\documentclass[twoside,leqno, twocolumn]{article}
\usepackage{cite}
\usepackage{amssymb}
\usepackage{graphicx}
\usepackage{url}
\usepackage{subcaption}
\usepackage[letterpaper]{geometry}
\usepackage{ltexpprt}
\usepackage{hyperref}
\usepackage{listings}
\usepackage{float}
\usepackage{wrapfig}
\usepackage{import}
\usepackage{booktabs}
\usepackage{siunitx}
\usepackage{algorithm}
\usepackage{algorithmic}
\usepackage{amsmath}
\usepackage{footmisc}
\usepackage{multirow}

\newcommand{\email}[1]{\protect\href{mailto:#1}{#1}}

\title{\Large Tuning Spectral Element Preconditioners for Parallel Scalability on GPUs}

\author{Malachi Phillips\thanks{Department of Computer Science, University of Illinois at Urbana-Champaign, Urbana IL 61801
  (\email{malachi2@illinois.edu}).
  }
\and Stefan Kerkemeier\thanks{Mathematics and Computer Science, Argonne National Laboratory, Lemont, IL 60439
  (\email{kerkemeier@anl.gov}).
  }
\and Paul Fischer$^{*\dagger}$\thanks{Department of Mechanical Science and Engineering,
University of Illinois at Urbana-Champaign, Urbana IL 61801
(\email{fischerp@illinois.edu}).
  }
}

\begin{document}
\date{}
\maketitle

\begin{abstract}
\small\baselineskip=9pt
The Poisson pressure solve resulting from the spectral element discretization
of the incompressible Navier-Stokes equation requires fast, robust, and
scalable preconditioning.  In the current work, a parallel scaling study of
Chebyshev-accelerated Schwarz and Jacobi preconditioning schemes is presented,
with special focus on GPU architectures, such as OLCF's Summit.  Convergence
properties of the Chebyshev-accelerated schemes are compared with alternative
methods, such as low-order preconditioners combined with algebraic multigrid.
Performance and scalability results are presented for a variety of
preconditioner and solver settings. The authors demonstrate that
Chebyshev-accelerated-Schwarz methods provide a robust and effective smoothing
strategy when using $p$-multigrid as a preconditioner in a Krylov-subspace
projector.
The variety of cases to be addressed, on a
wide range of processor counts, suggests that
performance can be enhanced by automated run-time selection of the
preconditioner and associated parameters.
\end{abstract}

\section{Introduction}
\label{introduction}
In fluid flow simulations, the incompressibility constraint is often used to
bypass fast acoustic waves and thereby allow the solution to evolve on a convective
time-scale that is most relevant for many engineering problems.  This model
leads to a Poisson problem for the pressure that is invariably the stiffest
substep in time-advancement of the Navier-Stokes equations. For large 3D
problems, which mandate the use of iterative methods, the pressure solve
thus typically encompasses the majority of the solution time.  

For the spectral element (SE) discretization, the Poisson system matrix
contains $\mathcal O(Ep^6)$ nonzeros for $E$ elements with a polynomial degree
of $p$ (i.e., approximately $n \approx Ep^3$ unknowns.)  Through the use of
tensor-product-sum factorization, however, the SE matrix-vector product can be
effected in only $\approx 7E(p+1)^3$ reads and $12E(p+1)^4$ operations, even in
the case of complex geometries \cite{parter_spectral_1979,deville_high-order_2002}.  
Consequently, the key to fast SE-based flow simulations is to find effective
preconditioners tailored to this discretization.  Many methods, including
geometric $p$-multigrid approaches with pointwise Jacobi and
Chebyshev-accelerated Jacobi smoothers
\cite{sundar_comparison_2015,adams_parallel_2003,Kronbichler2019}, 
geometric $p$-multigrid with overlapping Schwarz smoothers
\cite{lottes_hybrid_2005,stiller_nonuniformly_2017,loisel_hybrid_2008}, and
preconditioning via low-order discretizations
\cite{parter_spectral_1979,olson_algebraic_2007,bello-maldonado_scalable_2019,canuto_finite-element_2010},
have been considered.

In this work, we explore the parallel performance of these methods and extend
the Schwarz-smoothing based $p$-multigrid of \cite{lottes_hybrid_2005} to
support restrictive-additive Schwarz ideas of Cai and Sarkis \cite{ras99} as
well as Chebyshev-accelerated smoothing.  In addition, we extend the low-order
preconditioning strategy of \cite{bello-maldonado_scalable_2019} to run on GPU
architectures through the use of AmgX \cite{naumov_amgx_2015} as an algebraric
multigrid (AMG) solver for the sparse system.  Numerical results are shown for
the SE-based pressure Poisson problem as well as for the Navier-Stokes
equation.  All methods considered are implemented by the authors in the
scalable open-source CFD code, nekRS \cite{fischer_nekrs_2021}.  
nekRS started as a fork of libParanumal \cite{libp} and uses
highly optimized kernels based on the Open Concurrent Compute Abstraction (OCCA) \cite{medina2014occa}.
Special focus is given to performance and scalability on large-scale GPU-based
platforms such as OLCF's Summit.

The structure of this paper is as follows.  Section \ref{implementation}
outlines the spectral element formulation for the Poisson problem in $\mathbb
R^3$, as well as describe the various preconditioners considered.  A brief
description of several model probems of interest are presented in section
\ref{cases}.  Numerical results are highlighted in section \ref{results}.
Finally, a brief summary of the survey of solver techniques considered is
provided in section \ref{conclusions}.

\section{Background and Implementation}
\label{implementation}
\def\dO{\partial \Omega}
\def\Oh{{\hat \Omega}}
\def\uu{{\underline u}}
\def\ur{{\underline r}}
\def\us{{\underline s}}
\def\uf{{\underline f}}
\def\tlam {{\tilde \lambda}}

We introduce here basic aspects of the SE Poisson discretization and associated
preconditioners.

\subsection{SE Poisson Discretization}

Consider the Poisson equation in $\mathbb R^3$,
\begin{align}
  -\nabla^2 u = f \text{ for } u,f \in \Omega \subset \mathbb R^3 \mapsto \mathbb R.
  \label{eqn:poisson}
\end{align}
Boundary conditions for the pressure Poisson equation are either periodic,
Neumann, or Dirichlet, with the latter typically applicable only on a small
subset of the domain boundary, $\dO_D$, correspondong to an outflow condition.
Consequently, Neumann conditions apply over the majority (or all) of the domain
boundary, which makes the pressure problem more challenging than the standard
all-Dirichlet case.

The SE discretization of (\ref{eqn:poisson}) is based on the weak form:
{\em  Find $u \in X_0^p$ such that,}
\begin{align}
  (\nabla v, \nabla u)_p = (v,f)_p 
\;\;
\forall v \in X_0^p,
\end{align}
where $X_0^p$ is a finite-dimensional approximation comprising the basis
functions used in the SE discretization, ${\phi_j(\mathbf x)}$, $j=1,\ldots,n$,
that vanish on $\dO_D$ and $(\cdot,\cdot)_p$ represents the discrete $L^2$
inner product based on Gauss-Lobatto-Legendre quadrature in the reference
element, $\Oh:=[-1,1]^3$. 
The basis functions allow us to represent the solution, $u$, as $u(\mathbf x) =
\sum u_j \phi_j(\mathbf x)$, leading to a linear system of unknown basis
coefficients,
\begin{align}
  A \underline u = B \underline f,
\end{align}
with respective mass- and stiffness-matrix entries, 
$B_{ij} := (\phi_i,\phi_j)_p$, and 
$A_{ij} := (\nabla \phi_i, \nabla \phi_j)_p$.  

$\Omega$ is tesellated into nonoverlapping hexahedral elements, $\Omega^e$,
for $e=1,\ldots,E$, with isoparametric mappings from $\Oh$ to $\Omega^e$
provided by
$\mathbf x^e(r,s,t) = \sum_{i,j,k}\mathbf x_{ijk}^e h_i(r) h_j(s) h_k(t)$,
for $i,j,k\in[0,p]$.
Each $h_*(\xi)$ is a $p$th-order Lagrange cardinal polynomial on
the Gauss-Lobatto-Legendre (GLL) qudrature points, $\xi_j \in [-1,1]$.
Similarly, the test and trial functions $u,v$ are written in local form as
$u^e(r,s,t) = \sum_{i,j,k}u_{ijk}^e h_i(r) h_j(s) h_k(t)$.
Continuity is ensured across the interface between adjacent elements by
enforcing $u_{ijk}^e = u_{\hat i\hat j\hat k}^e$
when $\mathbf x_{ijk}^e = \mathbf x_{\hat i\hat j\hat k}^e$.
From this, a global-to-local degree-of-freedom mapping, 
$\uu := \{u_l \} \longrightarrow \uu_L := \{ u_{ijk}^e \}$,
can be represented by a Boolean matrix $Q$, such that $\uu_L = Q\uu$.
The assembled stiffness matrix is then $A=Q^TA_L Q$, where $A_L =
\operatorname{block-diag}(A^e)$ comprises the local stiffness matrices, $A^e$.
Similarly, $B=Q^TB_LQ$.  The SE formulation uses coincident GLL quadrature and
nodal points, such that $B^e$ is diagonal.  Moreover, $A^e$ is never formed, as
it would contain $O(p^6)$ nonzeros in the general case.  Rather, the
tensor-product-sum factorization \cite{parter_spectral_1979} allows for
$A\underline{u}$ to be evaluated in $O(Ep^4)$ time with $O(Ep^3)$ storage,
as described in detail in \cite{deville_high-order_2002}.

\subsection{Preconditioners}

In the current study, all preconditioners are applied in the context of
restarted GMRES.  Although $A$ is symmetric positive definite (SPD), many of the
preconditioners are asymmetric.  Further, a recent study 
\cite{fischer_highly_2021} has shown the benefits of projection-based GMRES
over flexible conjugate gradients (FCG) because the effectiveness of the overall
pressure solution strategy (which includes projection onto prior solutions
\cite{fischer_projection_1998}) generally ensures a low enough iteration
count, $k$, such that the $O(k^2)$ costs in GMRES are not overly onerous.
We note, however, that FCG may yield very low iteration counts
in certain cases (e.g., 1--2 iterations per step, as in \cite{fischer_nekrs_2021}),
in which case we use FCG rather than GMRES if it yields faster runtimes.

\subsubsection{SEMFEM} \label{sub:semfem}

In \cite{parter_spectral_1979}, Orszag suggested that constructing a sparse
preconditioner based on the low-order discretizations with nodes coinciding
with those of the high-order discretization would yield bounded condition
numbers and, under certain constraints, can yield $\kappa(M^{-1}A) \sim \pi^2 /
4$ for second-order Dirichlet problems.  This observation has led to the
development of preconditioning techniques based on solving the resulting low-order system
\cite{olson_algebraic_2007,bello-maldonado_scalable_2019,canuto_finite-element_2010}.

In the current work, we employ the same low-order discretization
considered in \cite{bello-maldonado_scalable_2019}.  Each of the vertices of
the hexahedral element is used to form one low-order, tetrahedral element,
resulting in a total of eight low-order elements for each GLL sub-volume
in each of the high-order hexahedral elements.
This low-order discretization is then used to form the sparse operator, $A_F$.
The so-called weak preconditioner, $A_F^{-1}$, is used to precondition the system.
Algebraric multigrid (AMG), implemented in CUDA in AmgX \cite{naumov_amgx_2015},
is used with the following setup to solve the low order system: \\[-2.5ex]
\begin{itemize}
  \item PMIS coarsening                              \\[-3.90ex]
  \item 0.25 strength threshold                      \\[-3.90ex]
  \item Extended + i interpolation ($p_{max}=4$)     \\[-3.90ex]
  \item Damped Jacobi relaxation ($0.9$)             \\[-3.90ex]
  \item One V-cycle for preconditioning              \\[-3.90ex]
  \item Smoothing on the coarsest level              \\[-3.90ex]
\end{itemize}
We denote this preconditioning strategy as SEMFEM.

\subsubsection{$p$-multigrid, Schwarz Smoothers}

Another preconditioning strategy for the SE-based Poisson problem is to use
geometric $p$-multigrid (pMG).  The classical single pass V-cycle is
summarized in Algorithm \ref{alg:V-cycle}.  In our application, we
limit pMG preconditioning to a single V-cycle pass.  Alternative
strategies, such as F- or W-cycle multigrid are not considered.

\begin{algorithm}\small
   \caption{Single pass multigrid V-cycle}
   \label{alg:V-cycle}
   \begin{algorithmic}
    \STATE{$\mathbf x = \mathbf x + \mathbf M (\mathbf b - \mathbf A\mathbf x)$ // smooth}
    \STATE{$\mathbf r = \mathbf b - \mathbf A \mathbf x$ // re-evaluate residual}
    \STATE{$\mathbf r_C = \mathbf P^T \mathbf r$ // coarsen }
    \STATE{$\mathbf e_C = \mathbf A_C^{-1} \mathbf r_C$ // solve/re-apply V-cycle }
    \STATE{$\mathbf e = \mathbf P \mathbf e_C$ // prolongate }
    \STATE{$\mathbf x = \mathbf x + \mathbf e$ // update solution }
    \STATE{$\mathbf x = \mathbf x + \mathbf M(\mathbf b - A\mathbf x)$ // post smoothing }
   \end{algorithmic}
\end{algorithm}

\begin{figure} \centering
{\setlength{\unitlength}{1.000in}
\begin{picture}(3.2,1.5)(0,1)
   \put(0.00,1.3){\includegraphics[width=1.50in]{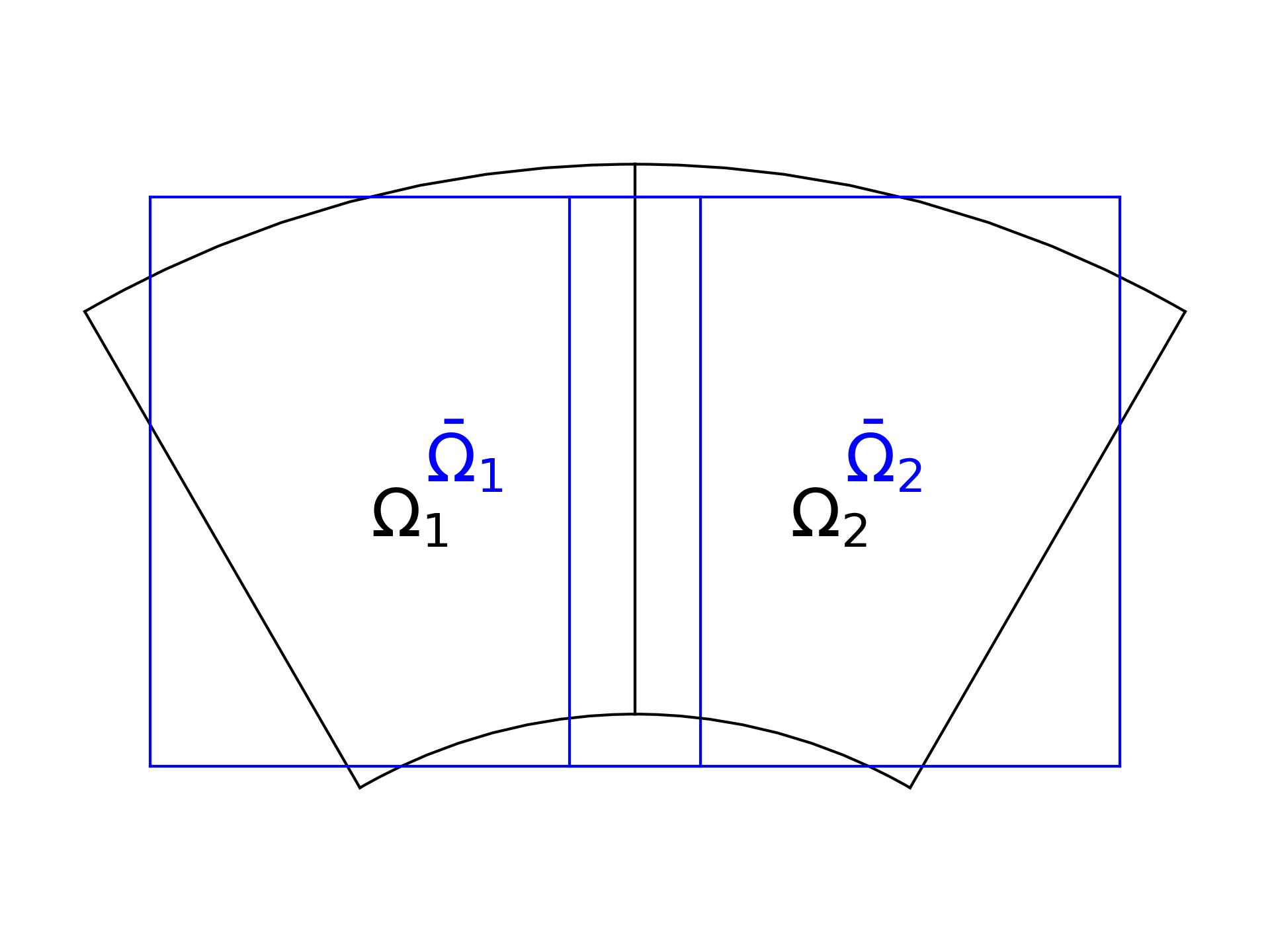}}
   \put(1.59,1.30){\includegraphics[width=1.50in]{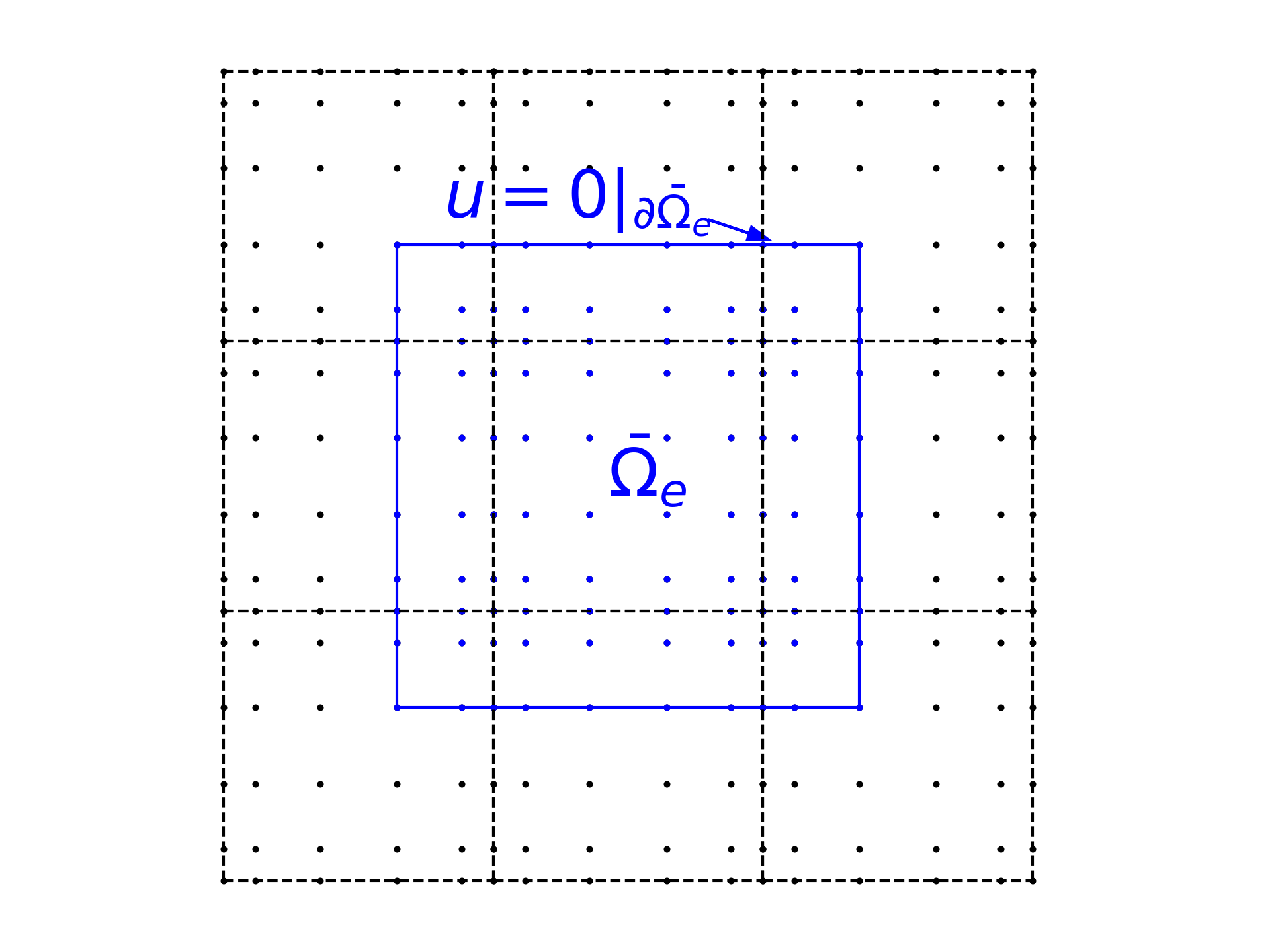}}

   \put(0.00,1.30){\large (a)}
   \put(1.55,1.33){\large (b)}

\end{picture}}
\vspace{-0.5cm}
\caption{
\small
(a) Box-like approximation (b) overlapping domain
\label{fig:schwarz-overlapping}
}
\end{figure}

The SE-based additive Schwarz method (ASM) presented in
\cite{lottes_hybrid_2005,loisel_hybrid_2008} solves local Poisson
problems on subdomains that are extensions of the spectral elements.
The formal definition of the ASM preconditioner (or pMG smoother) is
\begin{eqnarray} \label{eq:asm}
  \us &=& \sum_{e=1}^E W_e R_e^T {\bar A}_e^{-1} R_e \ur,
\end{eqnarray}
where $R_e$ is the restriction matrix that extracts nodal values
of the residual vector that correspond to each overlapping domain,
as indicated in Fig. \ref{fig:schwarz-overlapping}b.
To improve the smoothing properties of the ASM, we introduce
the diagonal weight matrix, $W_e$, which scales each nodal value by
the inverse of the number of subdomains that share that node.
Although it compromises symmetry, post-multiplication by $W_e$ was found to yield
superior results to pre- and post-multiplication by $W_e^{\frac{1}{2}}$
\cite{stiller_nonuniformly_2017,lottes_hybrid_2005}.

In a standard Galerkin ASM formulation, one would use ${\bar A}_{e} = R_e A
R_e^T$, but such an approach would compromise the $O(p^3)$ storage complexity
of the SE method.  To construct fast inverses for ${\bar A}_e$, we approximate
each deformed element as a simple box-like geometry, as demonstrated in Fig.
\ref{fig:schwarz-overlapping}a.  These boxes are then extended by a single 
degree-of-freedom in each spatial dimension to form overlapping subdomains with
${\bar p}^3=(p+3)^3$ interior degrees-of-freedom in each domain.
(See Fig.  \ref{fig:schwarz-overlapping}).  
The approximate box domain enables the use of the fast diagonalization method
(FDM) to solve for each of the overlapping subdomains, which can be applied in
$O(Ep^4)$ time in $\mathbb R^3$\cite{lottes_hybrid_2005}.  The
extended-box Poisson operator is separable, with the 3D form
\begin{equation}\label{eq:tensor-prod-poisson}
{\bar A} = 
B_z \otimes B_y \otimes A_x+B_z \otimes A_y \otimes B_x+A_z \otimes B_y \otimes B_x,
\end{equation}
where each $B_*,A_*$ represents the extended 1D mass-stiffness matrix pairs
along the given dimension \cite{lottes_hybrid_2005}.  The FDM begins with a
preprocessing step of solving a series of small,
${\bar p}\times {\bar p}$, generalized eigenvalue problems, 
\begin{equation}\label{eq:generalized-eig}
A_* \underline{s}_i = \lambda_i B_* \underline{s}_i
\end{equation}
and defining
$S_*=(\underline{s}_1\ldots\underline{s}_{\bar p})$ and
$\Lambda_*=\operatorname{diag}(\lambda_i)$, to yield the
similarity transforms
\begin{equation}\label{eq:similarity}
S_*^T A_* S_* = \Lambda_*, \;\; S_*^T B_* S_* = I.
\end{equation}
From these, the inverse of the local Schwarz operator is \\[-2.5ex]
\begin{equation}\label{eq:inverse-tensor-prod-poisson}
{\bar A}^{-1} = 
(S_z\otimes S_y\otimes S_x) D^{-1} (S_z^T \otimes S_y^T \otimes S_x^T),
\end{equation}
where
\begin{equation}
D=I\otimes I\otimes \Lambda_x+I\otimes \Lambda_y \otimes I+\Lambda_z \otimes I \otimes I
\end{equation}
is a diagonal matrix.  This process is repeated for each element, at each 
multigrid level save for the coarsest one.  Note that the
per-element storage is only $3{\bar p}^2$ for the $S_*$ matrices and ${\bar
p}^3$ for $D$.  At each multigrid level, the local subdomain solves are used as
a smoother.  On the coarsest level ($p=1$), however, BoomerAMG
\cite{henson_boomeramg_2002} is used to solve the system with the same
parameters as the AmgX solver in section \ref{sub:semfem}, except using
Chebyshev smoothing.  Unless otherwise noted, all pMG preconditioners use a
single BoomerAMG V-cycle iteration in the coarse-grid solve.

Presently, we also consider a restrictive additive Schwarz (RAS) version
of (\ref{eq:asm}), wherein overlapping values are not added after the action 
of the local FDM solve, following \cite{ras99}.  RAS has the added benefit of
reducing the amount of communication required in the smoother. (Formally,
RAS can be implemented by simply changing $W_e$.)
Note that, in the case of ASM and RAS smoothers without Chebyshev-acceleration,
an additive V-cycle with no post-smoothing is used to avoid residual 
re-evaluation in Algorithm \ref{alg:V-cycle} \cite{fischer04}.

\subsubsection{Chebyshev Acceleration} \label{sub:cheby}
A notable improvement over standard Jacobi-smoothed multigrid is to use 
Chebyshev-acceleration \cite{sundar_comparison_2015,adams_parallel_2003},
as described in Algorithm \ref{alg:cheby} for a given surrogate smoother, $S$. 
While $S$ is typically based on Jacobi smoothing (e.g.,\cite{Kronbichler2019}),
it is also possible to consider the use of overlapping-Schwarz (\ref{eq:asm})
as the smoother, which is a new approach that we explore here.

\begin{algorithm}\small
  \caption{Chebyshev smoother}
  \label{alg:cheby}
  \begin{algorithmic}
    \STATE{$\theta = \dfrac 1 2 (\lambda_{max}+\lambda_{min})$, $\delta= \dfrac 1 2 (\lambda_{max}-\lambda_{min})$, $\sigma = \dfrac \theta \delta$, $\rho_1 = \dfrac 1 \sigma$}
    \STATE{$ r =  S( b -  A  x)$, $ d_1 = \dfrac 1 \theta r$, $ x_1 =  0$}
    \FOR{$k=1,\dots,\texttt{chebyshevOrder}$}
      \STATE{$ x_{k+1} =  x_k +  d_k$}
      \STATE{$ r_{k+1} =  r_k -  S  A  d_k$}
      \STATE{$\rho_{k+1} = \dfrac{1}{2\sigma-\rho_k}$}
      \STATE{$ d_{k+1}=\rho_{k+1}\rho_k  d_k + \dfrac {2\rho_{k+1}}{\delta}  r_{k+1}$}
    \ENDFOR
    \STATE{$ x_{k+1} =  x_k +  d_k$}
    \RETURN{$ x_{k+1}$}
  \end{algorithmic}
\end{algorithm}

Algorithm \ref{alg:cheby} requires approximate spectral bounds, 
$(\lambda_{min}, \lambda_{max})$, of the smoothed operator, $SA$.  
Using 10 rounds of Arnoldi iteration, we generate $\tilde\lambda$ as an initial
proxy for $\lambda_{max}$.   The bounds employed in Algorithm \ref{alg:cheby} 
are then determined through sensitivity analysis similar to that performed by
Adams and coworkers \cite{adams_parallel_2003}.  The analysis is conducted
using a challenging model problem, namely the Kershaw case of subsection
\ref{kershaw}, with $\varepsilon=0.3$, $E=24^3,p=7$,
and relative ($2$-norm) residual tolerance of $10^{-8}$.
A second-order Chebyshev-accelerated ASM smoother is used in the $p=7,3,1$
pMG V-cycle preconditioner, which we denote as the current default
preconditioner in nekRS.

The results of the sensitivity study are shown in 
Table \ref{table:kershaw-eig-multiplier-0.3}.
The required iteration counts are seen to be sensitive to underestimation of
$\lambda_{max}$, as also found by Adams and coworkers \cite{adams_parallel_2003}.
On the other hand, results are less sensitive to the minimum eigenvalue
estimate, provided $\lambda_{min} > 0$, as this delimits between the low
frequencies that are handled by the coarse grid correction and the high
frequencies that are eliminated by the action of the smoother.  
The choice of the bounds (relative to $\tlam$) is not universal:
(1/30,1.1) \cite{adams_parallel_2003},
(0.3,1) \cite{baker_multigrid_2011},
(0.25,1) \cite{sundar_comparison_2015},
and
(1/6,1) \cite{zhukov_multigrid_2015} have all been considered.
Based on the results of Table \ref{table:kershaw-eig-multiplier-0.3}, we set
$(\lambda_{min},\lambda_{max})=(0.1,1.1)\tlam$ as a conservative
choice.

\def\tlam {{\tilde \lambda}}
%\vspace*{.01in}
%\def\tlam {{\vspace*{.11in}{\tilde \lambda}}}
\begin{table}
\small
\centering
\begin{tabular}{||c|c c c c c c||}
  \hline
  $\lambda_{min}\backslash \lambda_{max}$ & \rule{0pt}{2.5ex}
$0.9 \tlam$ & $0.95 \tlam$ & $1.0 \tlam$ & $1.1 \tlam$ & $1.2 \tlam$ & $1.3 \tlam$\\ 
\hline\hline \rule{0pt}{2.5ex}
$0.0   \tlam$ & - & - & - & - & - & - \\
$0.025 \tlam$ & - & - & - & 110  & 64   & 48 \\
$0.05  \tlam$ & - & - & - & 50   & 40   & 38 \\
$0.1   \tlam$ & - & - & 124  & 40   & 38   & 38 \\
$0.2   \tlam$ & 159  & 45   & 43   & 42   & 43   & 44 \\
$0.25  \tlam$ & 47   & 45   & 44   & 44   & 45   & 46 \\
\hline
\end{tabular}
\caption{\small
%Kershaw ($\varepsilon=0.3, E=24^3, p=7$) iteration counts required to reduce
%the residual by $10^8$.  Preconditioner is $p=7,3,1$ p-multigrid V-cycle with a
%2nd-order Chebyshev-accelerated ASM smoother.
Iteration counts for the ($\lambda_{min},\lambda_{max}$) sensitivity study.
Omitted entries failed to converge in 1000 iterations.
%Smoothing is based on 2nd-order Chebyshev-accelerated ASM with the pMG schedule
%$p=7,3,1$.  The test case is the Kershaw problem ($\varepsilon=0.3, E=24^3,
%p =7$) with relative residual tolerance $10^{-8}$.
\label{table:kershaw-eig-multiplier-0.3}}
\end{table}

\section{Test Cases}
\label{cases}
\def\dt{ \Delta t }

We describe four model problems that are used to test the SE preconditioners.
The first is a stand-alone Poisson solve, using variations of the Kershaw mesh.
The others are modest-scale Navier-Stokes problems, where the pressure Poisson
problem is solved over multiple timesteps.  The problem sizes are listed in
Table \ref{table:problem-sizes} and range from relatively small ($n$=21M points) to moderately large
($n$=645M).\footnote{Larger cases for recent full-scale runs on Summit with 
$n$=51B are reported in \cite{fischer_highly_2021}.}

\subsection{Poisson} \label{kershaw}

\begin{figure}
  \begin{subfigure}{0.15\textwidth}
    \includegraphics[width=\textwidth]{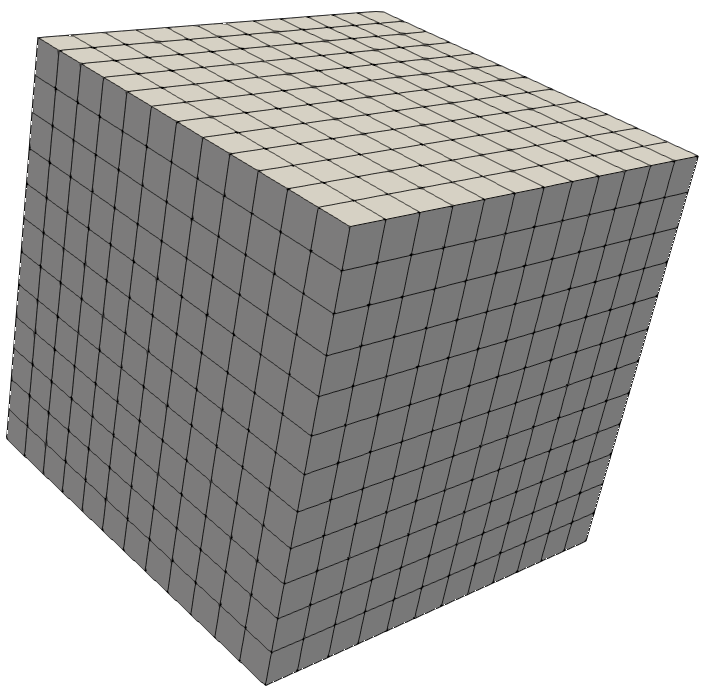}
  \end{subfigure}
  \hfill
  \begin{subfigure}{0.15\textwidth}
    \includegraphics[width=\textwidth]{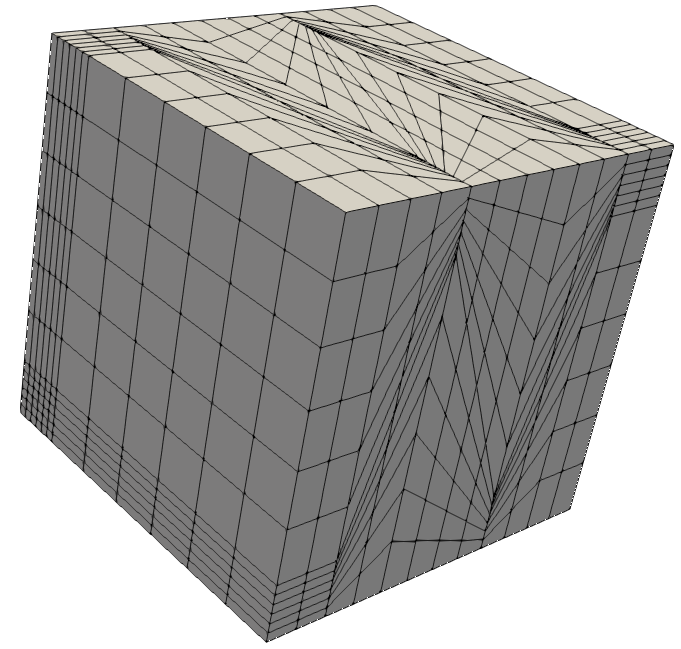}
  \end{subfigure}
  \hfill
  \begin{subfigure}{0.15\textwidth}
    \includegraphics[width=\textwidth]{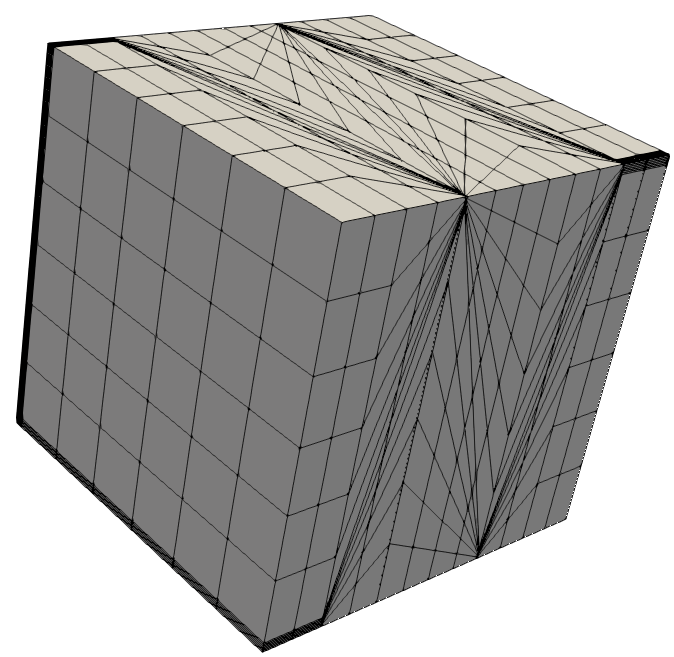}
  \end{subfigure}
  \caption{
    \label{fig:kershaws}
    Kershaw, $E=12^3,p=1$. $\varepsilon=1.0,0.3,0.05$.
  }
\end{figure}

The Kershaw family of meshes \cite{kolev_ceed_2021,kershaw_differencing_1981}.
has been proposed as the basis for a high-order Poisson-solver benchmark by
Center for Efficient Exascale Discretization (CEED) within the DOE Exascale
Computing Project (ECP).  This family is parametrized by an anisotropy measure,
$\varepsilon =\varepsilon_y=\varepsilon_z \in (0,1]$, that determines the
degree of deformation in the $y$ and $z$ directions.  As $\varepsilon$
decreases, the mesh deformation and aspect ratio increase along with it.  
The Kershaw mesh is shown in Fig. \ref{fig:kershaws} for several values
of $\varepsilon$.  The domain $\Omega = [-1/2,1/2]^3$ with Dirichlet boundary
conditions on $\partial\Omega$.  The right hand side for (\ref{eqn:poisson}) is
set to
\begin{equation}\label{eq:kershaw-rhs}
  f(x,y,z) = 3\pi^2 \sin{(\pi x)}\sin{(\pi y)}\sin{(\pi z)} + g,%(x,y,z),
\end{equation}
where $g(x,y,z)$ is a random, continous vector vanishing on $\partial\Omega$.
The linear solver terminates after reaching a relative residual reduction of
$10^{-8}$.  Since this test case solves the Poisson equation, there is no
timestepper needed for the model problem.

\subsection{Navier-Stokes} \label{navier-stokes}

\begin{table}
\centering
\begin{tabular}{||c| c c c ||}
  \hline
  Case Name & $E$ & $p$ & $n$\\
  \hline\hline
  146 pebble  (Fig. \ref{fig:ns_cases}a) & 62K & 7 & 21M\\
  1568 pebble (Fig. \ref{fig:ns_cases}b) & 524K & 7 & 180M\\
  67 pebble   (Fig. \ref{fig:ns_cases}c) & 122K & 7 & 42M\\
  Speed bump  (Fig. \ref{fig:ns_cases}d) & 885K & 9 & 645M\\
\hline
\end{tabular}
\caption{
  \small
  Problem discretization parameters.
  \label{table:problem-sizes}
}
\end{table}

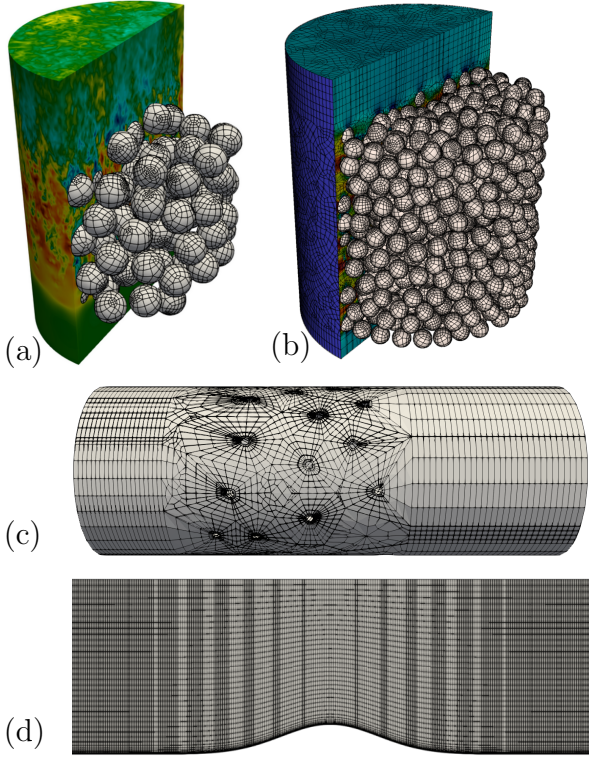
\begin{figure}[h] \centering
{\setlength{\unitlength}{0.9in}
\begin{picture}(3.2,4.4)(0,0)
   \put(0.00,2.27){\includegraphics[width=1.26in]{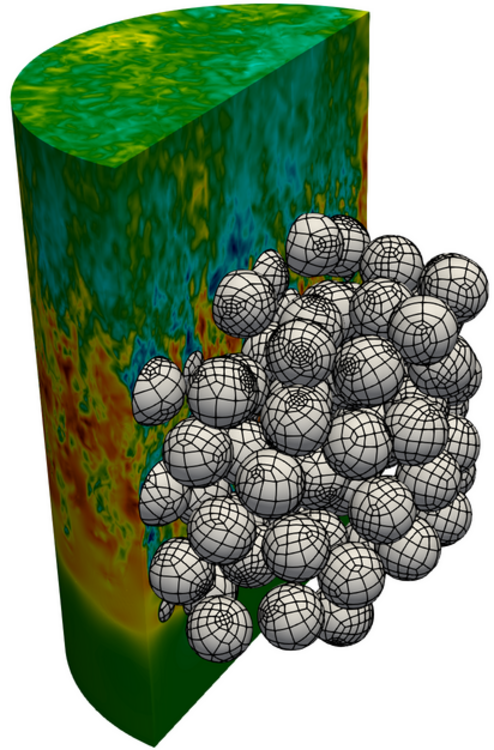}}
   \put(1.59,2.30){\includegraphics[width=1.44in]{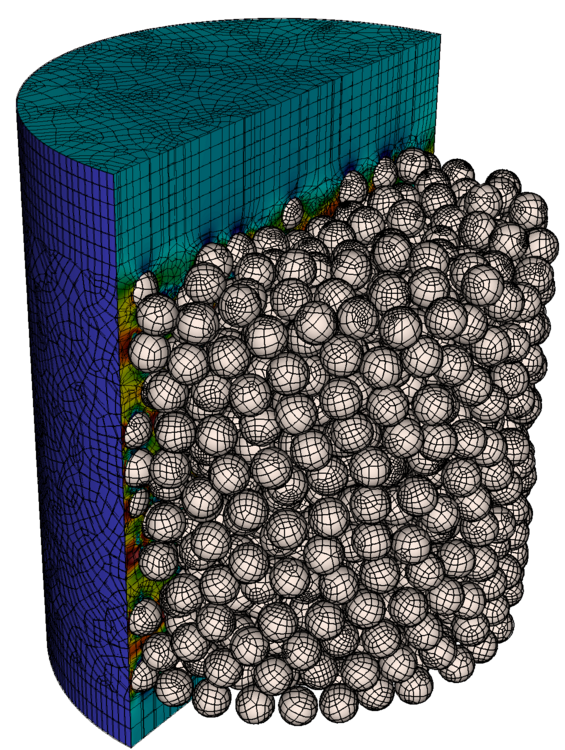}}
   \put(0.40,1.15){\includegraphics[width=2.70in]{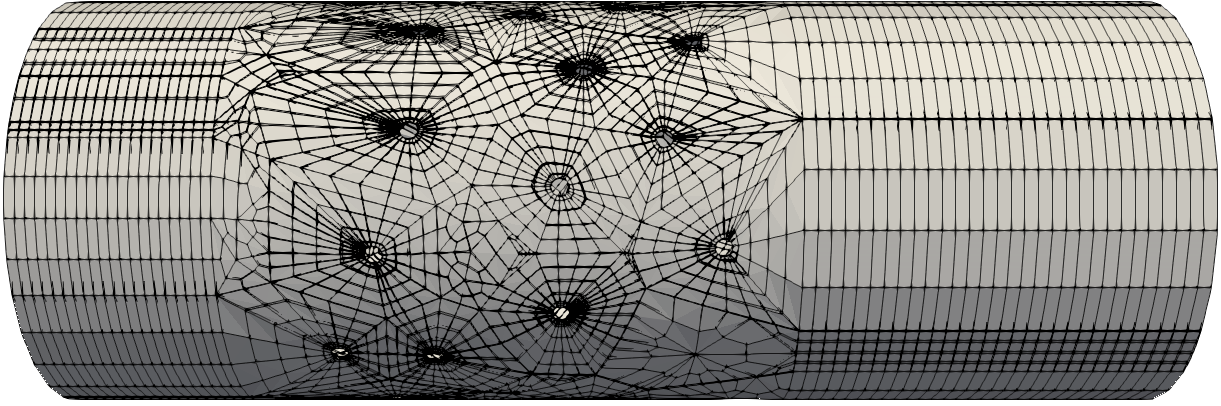}}
   \put(0.40,0.00){\includegraphics[width=2.70in]{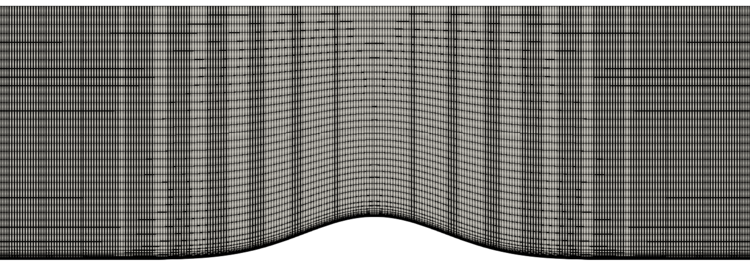}}

   \put(0.00,2.30){\large (a)}
   \put(1.55,2.33){\large (b)}
   \put(0.00,1.23){\large (c)}
   \put(0.00,0.08){\large (d)}

\end{picture}}
\caption{Navier-Stokes cases: pebble-beds with 
(a) 146, (b) 1568, and (c) 67 spheres;
(d) Boeing speed bump.
\label{fig:ns_cases}
}
\end{figure}

For the pressure-Poisson tests, four flow cases are considered, as depicted
in Fig.  \ref{fig:ns_cases}.  
The first three cases corresponds to turbulent flow through a cylindrical packed-bed
with 146, 1568, and 67 spherical pebbles.
The 146 and 1568 pebble cases are from Lan and coworkers \cite{lan_all-hex_2021}.
The 67 pebble case is constructed using a tet-to-hex meshing strategy by Yuan and coworkers \cite{yuan2020spectral}.
The first two bed flows are at Reynolds number $Re_D=5000$, based on sphere
diameter, $D$, while the 67 pebble case is at Reynolds number $Re_D=1460$.  Time advancement is based on a two-stage 2nd-order
characteristics timestepper with CFL=4 
($\dt = 2\times10^{-3}$
$\dt = 5\times 10^{-4}$, and
$\dt = 5\times 10^{-5}$ for the 146, 1568, and 67 pebble cases). 
An absolute pressure solver tolerance of $10^{-4}$ is used.  A restart at
$t=10$, $t=20$, and $t=10.6$ convective time units is used for the 146, 1568, and 67 pebble
cases, respectively, to provide an initially turbulent flow.   

The fourth case, shown in Fig. \ref{fig:ns_cases}d, is a direct numerical
simulation (DNS) of seperated turbulent flow over a speed bump at $Re = 10^6$.
This test case was designed by Boeing to provide a flow that exhibits
separation.  A DNS of the full 3D geometry, however, remains difficult
\cite{shur_direct_2021}.  Therefore, this smaller example proves a useful
application for benchmarking solver performance.  This case uses a 2nd-order
timestepper with CFL=0.8 ($\Delta t = 4.5\times 10^{-6}$) and an absolute
pressure-solve tolerance of $10^{-5}$.  A restart at $t=5.6$ convective time
units is used for the initial condition.

In all cases, solver results are collected over 2000 timesteps.  At each step,
the solution is projected onto a space of up to 10 prior solution vectors to
generate a high-quality initial guess, ${\bar \uu}$.  Projection is standard
practice in nekRS as it can reduce the initial residual by orders of magnitude
at the cost of just one or two matrix-vector products in $A$ per step
\cite{fischer_projection_1998}.

The perturbation solution, $\delta \uu := \uu - {\bar \uu}$, is
typically devoid of slowly evolving low wave-number content.  Moreover, the
initial residual is oftentimes sufficiently small that the solution converges
in $k < 5$ iterations, such that the $O(k^2)$ overhead of GMRES is small.
Testing the preconditioners under these conditions ensures that the conclusions
drawn are relevant to the application space.

\section{Results}
\label{results}
Here we consider the solver performance results for the test cases
of Section 3.  We assign a single MPI rank to each GPU
and denote the number of ranks as $P$. All runs are on Summit.
Each node on Summit consists of 42 IBM Power9 CPUs and 6 NVIDIA V100 GPUs.  
We use 6 GPUs per node unless $P < 6$.
In the following, we denote a pMG preconditioner using $\eta$-order Chebyshev-accelerated $\xi$ smoother with
a multigrid schedule of $\Pi$ as Cheby$\xi$($\eta$),$\Pi$.
A wide range of preconditioning strategies is considered.

\subsection{Kershaw Mesh}

The Kershaw study comprises six tests.  For each of two studies, we consider
the regular box case ($\varepsilon=1.0$), a moderately skewed case
($\varepsilon=0.3$), and a highly skewed case ($\varepsilon=0.05$).
The first study is a standard weak-scale test, where $P$ and
$E$ are increased, while the polynomial order is fixed at $p=7$ and the number
of gridpoints per GPU is set to $n/P = 2.67M$.  The range of processors is
$P$=6 to 384.  The second study is a test of the influence of polynomial order
on conditioning, with $P=24$ and $n/P=2.88M$ fixed, while $p$ ranges from 3 to 10.
Both cases use GMRES(20).  

\begin{figure*}
  \centering
  \includegraphics[width=\textwidth]{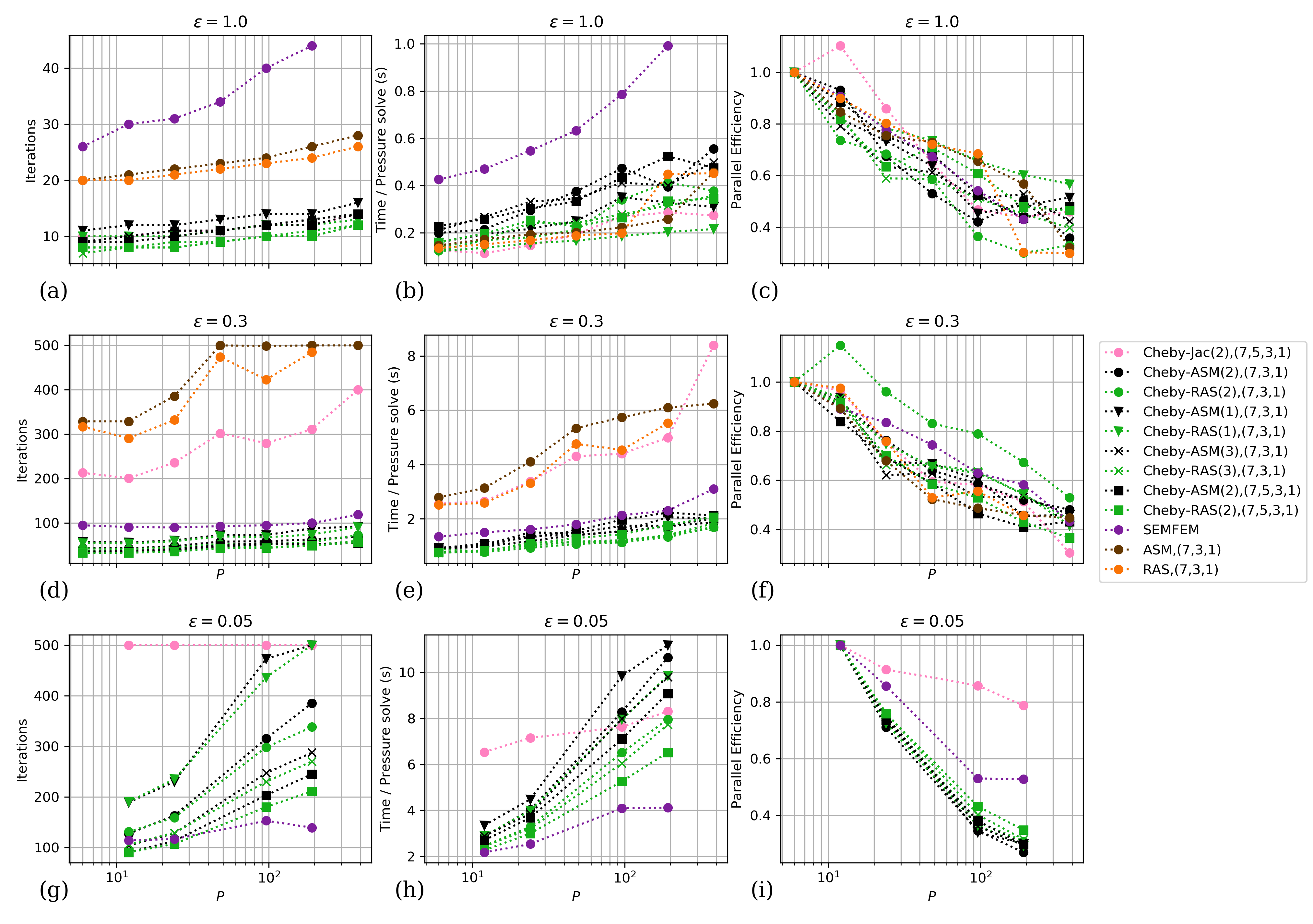}
  \vspace*{-1cm}
  \caption{
      \small
      Kershaw weak scaling, GMRES(20).
      $n/P = 2.67M$.
      \label{fig:kershaw-weak-scaling}
  }
\end{figure*}

The results of the weak-scaling study are shown in Fig.
\ref{fig:kershaw-weak-scaling}.  
For all values of $\varepsilon$, the iteration count exhibits a dependence 
on problem size, as seen in Fig.  \ref{fig:kershaw-weak-scaling}a,d,g,
especially in the highly skewed case ($\varepsilon = 0.05$).
The time-per-solve also increases with $n$, in part due to the increase
in iteration count, but also due to increased communication overhead as
$P$ increases. 
This trend is not necessarily monotonic, as shown in Fig. \ref{fig:kershaw-weak-scaling}c,f.
In the case of $\varepsilon=0.3$ Cheby-RAS(2),(7,3,1), a minor fluctuation in the iteration count
from the $P=6$ to $P=12$ case causes the greater than unity parallel efficiency.
For $\varepsilon=1.0$ Cheby-Jac(2),(7,5,3,1), the effects of system noise
on the $P=6$ run causes the greater than unity parallel efficiency
for the $P=12$ run.

Table \ref{table:max-neighbors} 
indicates the maximum number of of neighboring processors for the assembly
($QQ^T$) graph of $A$, which increases with $P$.  In addition, the number of
grid points per GPU ($n/P=2.67M$) is relatively low.  These two factors cause
an increased sensitivity of the problem to the additional communication
overhead as the number of GPUs is initially increased.  The number of
neighbors, however, will saturate at larger (i.e., production-level) processor
counts.

Lastly, the relative preconditioner performance depends on $\varepsilon$.
Fig. \ref{fig:kershaw-weak-scaling}b, e, h show that, for the easy $\varepsilon=1.0$ case,
a pMG scheme with a smoother that is cheap to apply is best, such as
Cheby-RAS(1),(7,3,1), Cheby-Jac(2),(7,5,3,1), and ASM,(7,3,1).
However, as $\varepsilon$ decreases, more robust pMG smoothers such as
Cheby-RAS(2),(7,3,1) and SEMFEM result in lower time to solution.
Once $\varepsilon = 0.05$ (Fig. \ref{fig:kershaw-weak-scaling}h), the problem is sufficiently challenging that
SEMFEM overtakes the pMG based preconditioning schemes.
This indicates that, in the highly skewed case in which the maximum element aspect ratio
increases, the pMG preconditioner is not as effective as the SEMFEM preconditioner.

\begin{figure*}
  \centering
  \includegraphics[width=\textwidth]{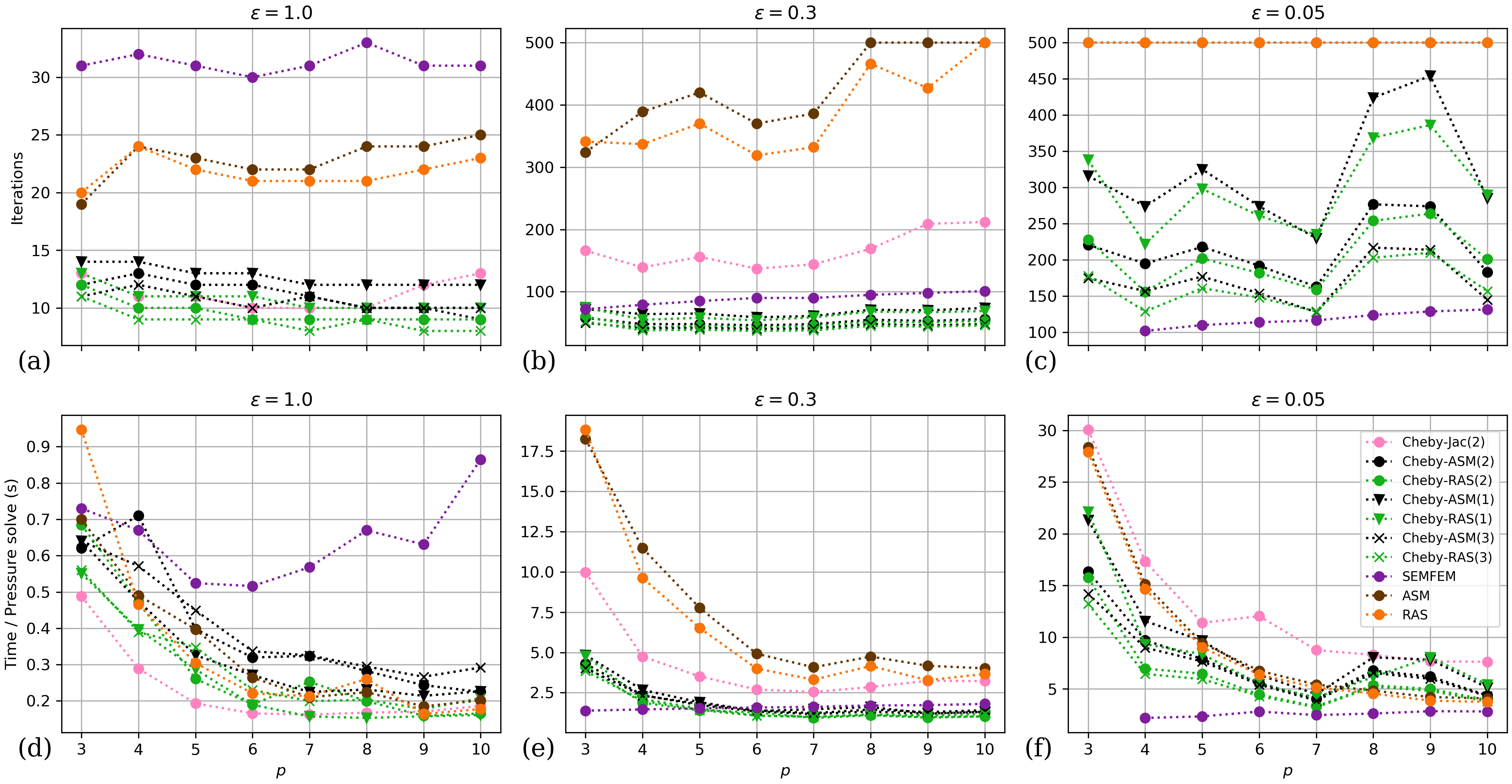}
  \vspace*{-0.8cm}
  \caption{
    \small
    Kershaw order dependence, GMRES(20).
    $n/P = 2.88M$.
    \label{fig:kershaw-order-dependence}
  }
\end{figure*}

\begin{table}
\centering
\begin{tabular}{|| c | c c c c c c c ||}
\hline
Nodes & 1 & 2 & 4 & 8 & 16 & 32 & 64\\
\hline
Max Neighbors & 5 & 11 & 16 & 24 & 20 & 24 & 29\\
\hline
\end{tabular}
%\vspace{-0.25cm}
\caption{
  \small
  Max neighbors, Kershaw.
  \label{table:max-neighbors}
}
\end{table}

The influence of polynomial order is illustrated in Fig.
\ref{fig:kershaw-order-dependence}.  For $\varepsilon=1.0$, iteration counts
are essentially $p$-independent, as seen in Fig.
\ref{fig:kershaw-order-dependence}a.  For $\varepsilon=0.3$, however, a slight
upward trend in the iteration count is observed for SEMFEM and for the pMG
preconditioners with ASM, RAS, and Chebyshev-Jacobi smoothing
(Fig.  \ref{fig:kershaw-order-dependence}b).
Similarly, $\varepsilon = 0.05$ exhibits a dependence between the iteration count and the polynomial order
for all preconditioners (Fig. \ref{fig:kershaw-order-dependence}c).
Figures \ref{fig:kershaw-order-dependence}d, e and f demonstrate that the {\em time} per
pressure solve is strongly dependent polynomial order.
For SEMFEM, there is an increase in time as a result of increased overhead for the underlying AMG solver.
For the pMG-based methods, higher orders are generally faster.  This
performance gain can be attributed to surface-to-volume effects in the
evaluation of the operators and smoothers.  Application of $A^e$ and (smoother)
$S^e$ is highly vectorizable, whereas application of $QQ^T$ (assembly) involves
a significant amount of indirect addressing.  The discrete surface to volume
ratio for a spectral element of order $p=7$ is 296/512 $\approx$ 60\%.  For
lower $p$ this value is larger.  This situation is exacerbated in the case of
Schwarz-based smoothers, where the overlap contributes substantially to the
work and communication for the small-element (i.e., low-$p$) cases.  In
addition, nekRS \cite{fischer_nekrs_2021} uses a single element per thread
block, which limits the amount of work available for a streaming multiprocessor
for relatively small polynomial orders.  

\subsection{Navier-Stokes Results}
\begin{figure}
  \centering
  \includegraphics[width=0.45\textwidth]{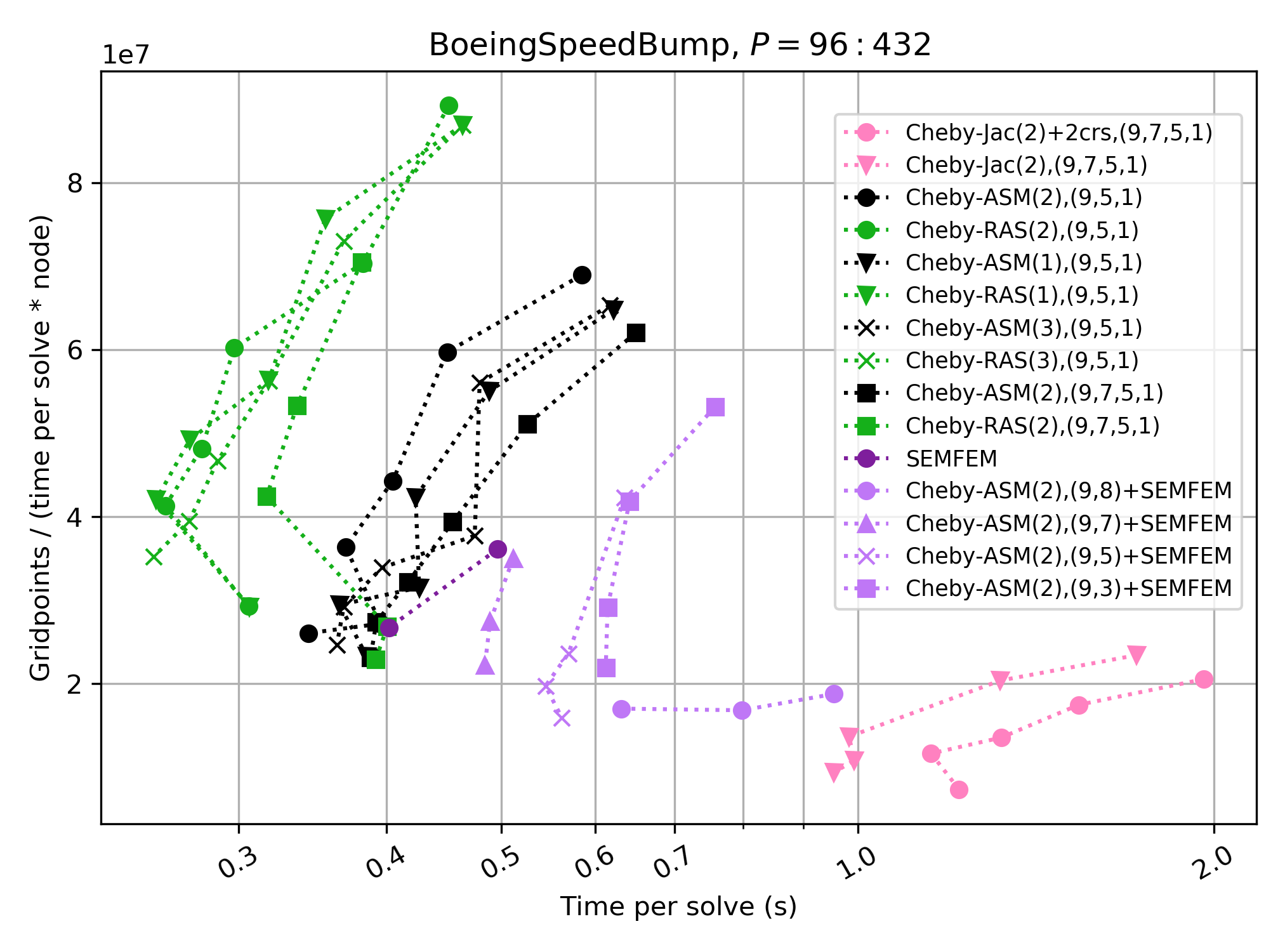}

  \vspace{-0.35cm}
  \caption{
    \small
    Strong scaling results on Summit for the Navier-Stokes case of Fig. \ref{fig:ns_cases}d.
  }
  \label{fig:scaling-study-bsb}

\end{figure}

\begin{figure*}
  \centering
  \includegraphics[width=\textwidth]{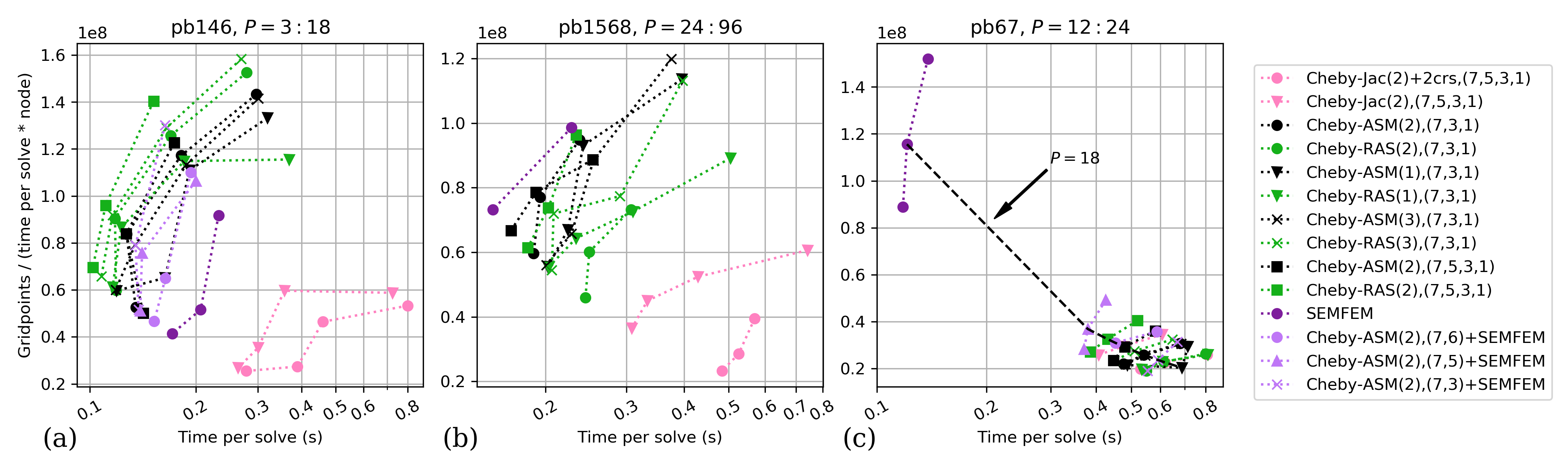}

  \vspace{-0.35cm}
  \caption{
    \small
    Strong scaling results on Summit for the Navier-Stokes cases of Fig. \ref{fig:ns_cases}a,b,c.
    Iso-processor count line illustrated in (c).
    A user running on a specified number of processors should use the lowest time-to-solution preconditioner
    along this line.
  }
  \label{fig:scaling-study-pb}

\end{figure*}

We consider scalability of nekRS for the cases of Fig. \ref{fig:ns_cases}.
All simulations except one use GMRES(15) with an initial guess generated by
$A$-conjugate projection onto 10 prior solutions \cite{fischer_projection_1998}.
Due to memory constraints, the 1568-pebble case with $P=24$ uses GMRES(10) with only
5 solution-projection vectors.
For each case, two pMG schedules are considered:
$(7,5,3,1)$ and $(7,3,1)$  for $p=7$;
and
$(9,7,5,1)$ and $(9,5,1)$ for $p=9$. 
Other parameters, such as the Chebyshev order and the number of
coarse grid BoomerAMG V-cycles are also varied.  
Results are shown in Figs. \ref{fig:scaling-study-bsb},\ref{fig:scaling-study-pb}.
The plots relate the effective work rate per node, measured as the
gridpoints $n$ (as shown in Table \ref{table:problem-sizes}) solved per second per node, to the time-to-solution.
The y-axis notes the drop in the relative work rate, which corresponds to a lower parallel efficiency, as the strong scale limit is reached,
while the x-axis denotes the time-to-solution.
Each node consists of 6 GPUs, hence $P=6 \times \text{nodes}$.

In all the performance tests conducted, the pMG preconditioner with
Chebyshev-Jacobi smoothing is outperformed by the other preconditioners,
whether using one or two V-cycle iterations in the AMG coarse-grid solve.  For
each case, the fastest preconditioner scheme varies.
In the 146 pebble case (Fig. \ref{fig:scaling-study-pb}a),
using Cheby-RAS(2),(7,5,3,1) yields the smallest time per pressure solve.
However, in the 1568 pebble case (Fig. \ref{fig:scaling-study-pb}b), SEMFEM is a moderate improvement over the
second best preconditioner, Cheby-ASM(2),(7,5,3,1).
pMG with a
$(9,5,1)$ schedule and Chebyshev-RAS (of any order) yield the best scalability
and lowest time per pressure solve for the Boeing speed bump case (Fig. \ref{fig:scaling-study-bsb}).
The Chebyshev-accelerated Schwarz schemes are not always the fastest, however.
For the 67 pebble case (Fig. \ref{fig:scaling-study-pb}c), Cheby-Jac(2),(7,5,3,1) is comparable to Cheby-RAS(2),(7,5,3,1)
and are the two fastest pMG based preconditioners.
However, SEMFEM is significantly faster than the other preconditioners for this case.

Also considered is a hyrbid two-level pMG/SEMFEM approach wherein SEMFEM is used as the solver for the coarse level.
For $p=7$, a $(7,6)$ schedule with 2nd order Chebyshev-accelerated ASM smoothing on the $p=7$ level
and SEMFEM solver on the $p=6$ level, denoted as Cheby-ASM(2),(7,6) + SEMFEM, is used.
Similarly, Cheb-ASM(2),(7,5) + SEMFEM and Cheb-ASM(2),(7,3) + SEMFEM are considered.
For $p=9$, a $(9,8)$, $(9,7)$, $(9,5)$, and $(9,3)$ hybrid pMG/SEMFEM approach is considered,
denoted as Cheb-ASM(2),(9,8) + SEMFEM, Cheb-ASM(2),(9,7) + SEMFEM, Cheb-ASM(2),(9,5) + SEMFEM, and Cheb-ASM(2),(9,3) + SEMFEM,
respectively.
In the pebble cases shown in Fig. \ref{fig:scaling-study-pb}a,c, this hybrid approach
performs somewhere between the SEMFEM and Cheby-ASM(2),(7,3,1) preconditioners.
In the Boeing speed bump case, however, Fig. \ref{fig:scaling-study-bsb} demonstrates
that this approach is not as performant as either the SEMFEM or Cheby-ASM(2),(9,5,1) preconditioners.

Solver performance degrades whenever $n/P$ is sufficiently small,
regardless of the solver considered.  For the various preconditioners
considered, at 2 nodes ($n/P=1.75M$), the strong-scale limit of 80\%
efficiency is far surpassed in the 146 pebble case.  This leaves the effective
strong scale limit at 1-2 nodes ($n/P=3.5$ to $1.75M$).
For the 1568 pebble case, 12 nodes ($n/P=2.5M$) yields a parallel
efficiency around 60-70\%, depending on the specific solver.
The 67 pebble cases reaches 70\% effieciency on 3 nodes ($n/P = 2.3M$).
The
parallel efficiency for the Boeing speed bump case for the fastest
time-to-solution preconditioners drops below 70\% when using
more than 48 nodes ($n/P=2.24M$).

\begin{figure}
  \centering
  \includegraphics[width=0.45\textwidth]{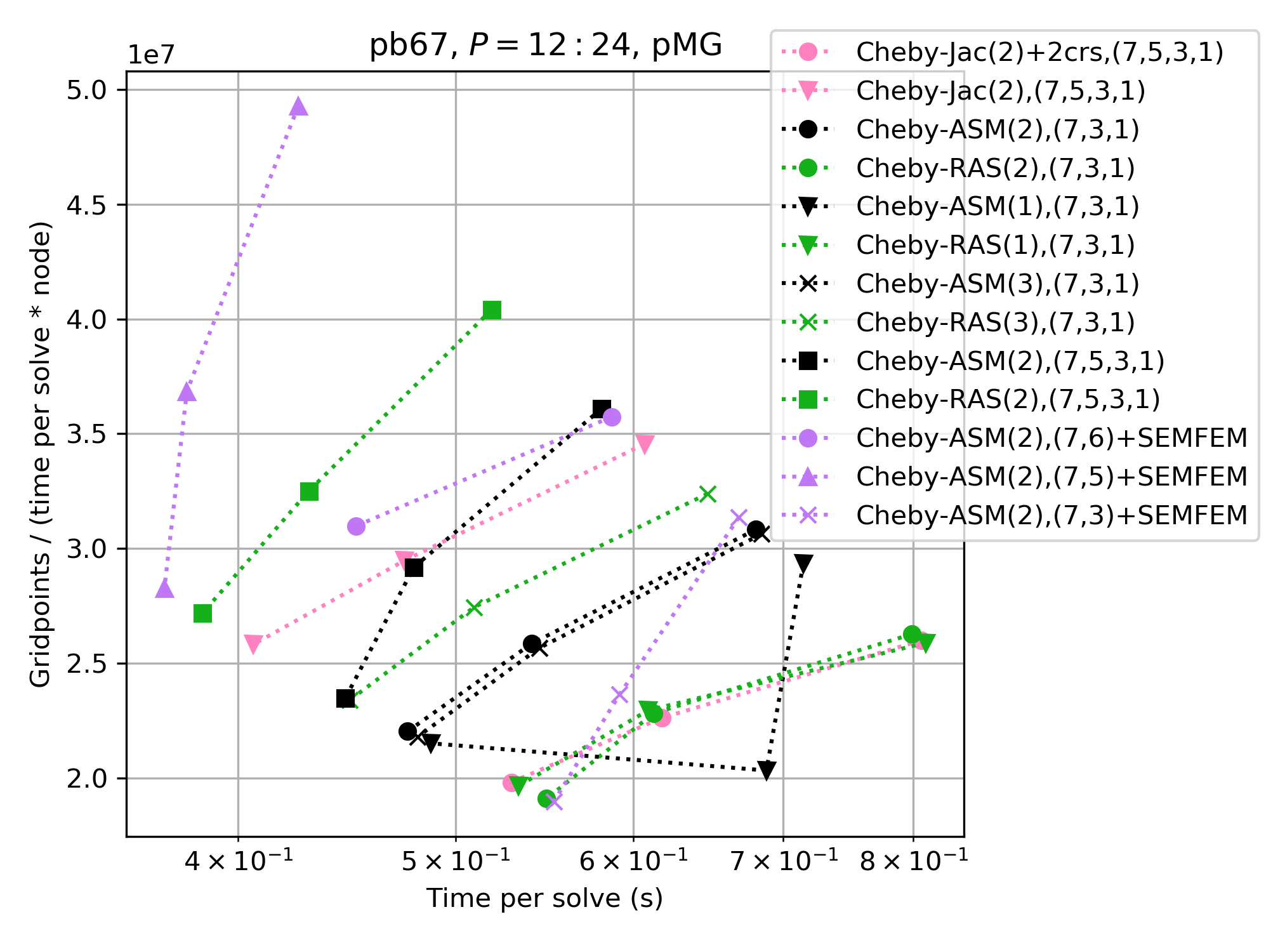}

  \vspace{-0.35cm}
  \caption{
    \small
    Same as Fig. \ref{fig:scaling-study-pb}c, pMG preconditioners only.
  }
  \label{fig:scaling-study-pb-67-zoom}

\end{figure}

While the effect of mesh quality metrics, such as the max aspect ratio and scaled Jacobian, on solver convergence has been studied by Mittal {\em et al.}
\cite{mittal2019mesh}, predicting the optimal preconditioner settings from
mesh quality metrics is not obvious.
While the maximum aspect ratio is an important metric for mesh quality, it alone cannot
explain the apparent poor performance of the pMG preconditioners in the 67 pebble case, Fig. \ref{fig:scaling-study-pb}c.
Consider, for example, that the Boeing speed bump case has a larger maximum aspect ratio (255)
than the 67 pebble case (204), but does not exhibit this poor performance in the pMG preconditioners.
The minimum scaled Jacobian, however, is at least an order of magnitude smaller in the 67 pebble case ($5.97\times 10 ^{-3}$) as compared
to the other cases (e.g., .996 for the Boeing speed bump case, $4.31\times 10^{-2}$ and $2.59\times 10^{-2}$ for the 146, 1568 pebble cases, respectively).
This may, in turn, explain why the pMG preconditioners in the 67 pebble case
were signficiantly suboptimal compared to the SEMFEM preconditioner.
This demonstrates, however, that a user cannot simply rely on, e.g., the maximum aspect ratio
when deciding whether or not to use SEMFEM as a preconditioner.
This inability to correctly identify preconditioner settings based on mesh quality metrics
alone motivates the introduction of an auto-tuner to choose preconditioner parameters
during runtime.
\begin{table*} \footnotesize
\centering
\begin{tabular}{||c| c c c ||}
  \hline
  %& \multicolumn{3}{c ||}{Mesh Quality Metrics}\\
  %\cline{2-4}
  Case Name & GLL Spacing (min/max) & Scaled Jacobian (min/max/avg) & Aspect Ratio (min/max/avg)\\
  %\hline\hline
  \hline \rule{0pt}{2.5ex}
  146 pebble
    &  $1.01\times 10^{-3}$ / .32
    & $4.31\times 10^{-2}$ / .977 / .419
    & 1.07 / 56.9 / 7.14 \\
  1568 pebble
    & $2.21\times 10^{-4}$ / .3
    & $2.59\times 10^{-2}$ / .99 / .371
    & 1.12 / 108 / 12.6 \\
  67 pebble
    & $4.02\times 10^{-5}$ / .145
    & $5.97\times 10^{-3}$ / .970 / .38
    & 1.17 / 204 / 13.2 \\
  Boeing Speed Bump
    & $8.34\times 10^{-7}$ / $2.99\times 10 ^{-3}$
    & .996 / 1 / 0.999
    & 6.25 / 255 / 28.1 \\
\hline
\end{tabular}
\vspace{-0.1cm}
\caption{
  \small
  Mesh quality metrics for cases from Fig. \ref{fig:ns_cases}.
  \label{table:mesh-metrics}
}
\end{table*}

%# bsb
%mesh metrics:
%GLL grid spacing min/max    : 8.34E-07 2.99E-03
%scaled Jacobian  min/max/avg: 9.96E-01 1.00E+00 9.99E-01
%aspect ratio     min/max/avg: 6.25E+00 2.55E+02 2.81E+01
%# pb67
%mesh metrics:
%GLL grid spacing min/max    : 4.02E-05 1.45E-01
%scaled Jacobian  min/max/avg: 5.97E-03 9.70E-01 3.80E-01
%aspect ratio     min/max/avg: 1.17E+00 2.04E+02 1.32E+01
%
%## pb146
%mesh metrics:
%GLL grid spacing min/max    : 1.01E-03 3.20E-01
%scaled Jacobian  min/max/avg: 4.31E-02 9.77E-01 4.19E-01
%aspect ratio     min/max/avg: 1.07E+00 5.69E+01 7.14E+00
%
%## pb1568
%mesh metrics:
%GLL grid spacing min/max    : 2.21E-04 3.02E-01
%scaled Jacobian  min/max/avg: 2.59E-02 9.91E-01 3.71E-01
%aspect ratio     min/max/avg: 1.12E+00 1.08E+02 1.26E+01
%

\section{Auto-tuning for Production Simulations}   \label{sec:online-preco} The
results of the preceding section provides a small window into the varieties of
performance behavior encountered in actual production cases, which span a large
range of problem sizes, domain topologies, mesh qualities.  Moreover,
production simulations are solved across a range of architectures having
varying on-node and network performance, interconnect topologies, and processor
counts.   One frequently encounters situations where certain communication
patterns might be slower under a particular MPI version on one platform versus
another.  Unless a developer has access to that platform, it's difficult to
measure and quantify the communication overhead.  Processor count alone can be
a major factor in preconditioner selection: large processor counts have
relatively high coarse-grid solve costs that can be mitigated by doing more
smoothing at the fine and intermediate levels.  How much more is the open
question.  The enormity of parameter space, particularly ``in the field''
(i.e., users working on unknown platforms) limits the effectiveness of standard
complexity analysis in selecting the optimal preconditioner for a given user's
application.

   From the user's perspective, there is only one application (at a time,
typically), and one processor count of interest.  Being able to provide
optimized performance---tuned to the application at hand, which includes the
processor count---is thus of paramount importance.  Auto-tuning provides an
effective way to deliver this performance.  Auto-tuning of preconditioners has
been considered in early work by Imakura {\em at al.}
\cite{imakura_auto-tuning_2012} and more recently by Yamada {\em at al.}
\cite{yamada2018preconditioner} and by Brown {\em et al.} \cite{brown_tuning_2021}.
The latter work couples their auto-tuning with local Fourier analysis to guide
the tuning process.

\begin{table*}[h] \footnotesize
\centering
\begin{tabular}{|l | l|}
\hline
& Parameters \\
\hline
Solver & preconditioned GMRES\\
\hline
Preconditioner & pMG, SEMFEM \\
\multirow{2}{*}{\hspace*{.15in} pMG} & $p=7$: Cheby-ASM(2),(7,3,1), Cheby-RAS(2), (7,3,1), Cheby-Jac(2), (7,5,3,1)\\
                                     & $p=9$: Cheby-ASM(2),(9,5,1), Cheby-RAS(2), (9,5,1), Cheby-Jac(2), (9,7,5,1)\\
\hspace*{0.3in} Coarse grid & single boomerAMG V-cycle\\
\hspace*{.15in} SEMFEM & single AmgX V-cycle\\
\hline
\end{tabular}

\vspace{-0.25cm}
\caption{
  \small
  Solver parameter space considered in the auto-tuner.
  A pMG preconditioner using an $\eta$-order Chebyshev-accelerated $\xi$ smoother with
  a multigrid schedule of $\Pi$ is denoted as Cheby-$\xi$($\eta$),$\Pi$.
  \label{table:auto-tuner-parameter-space}
}
\end{table*}

In large-scale fluid mechanics applications, auto-tuning overhead is typically
amortized over $10^4$--$10^5$ timesteps (i.e., pressure solves) {\em per run}
(and more, over an entire simulation campaign).  Moreover, auto-tuning is of
particular importance for problems at large processor counts because these
cases often have long queue times, which preclude making multiple job
submissions in order to tweak parameter settings.  Failure to optimize,
however, can result in significant opportunity costs.  For example, in
\cite{fischer_highly_2021} the authors realized a factor of 2.8 speedup in time-per-step for a
352,000-pebble-bed simulation ($n$=51B gridpoints) through a sequence of tuning
steps.  Even if there were 100 configurations in the preconditioner parameter
space, an auto-tuner could visit each of these in succession 5 times each within
the first 500 steps and have expended only a modest increment in overhead when
compared to the cost of submitting many jobs (for tuning) or to the cost of
10,000 steps for a production run.

As a preliminary step, the authors consider constructing
a small subset of the true search space (which could include, e.g., other cycles, schedules,
or smoothers) for use in a nascent auto-tuner.
Across all cases, with exception to the 67 pebble case, pMG preconditioning
with Chebyshev-accelerated ASM or RAS smoothing is the fastest solver or is comparable to SEMFEM.
The choice of Chebyshev order and multigrid schedule, moreover,
contributes only a modest $\approx$ 10-20\% improvement to the overall time-to-solution in most cases, all else being equal.
This makes the default Cheby-ASM(2),(7,3,1) or Cheby-ASM(2),(9,5,1) preconditioner
reasonably performant.
However, in order to avoid the situation encountered in the 67 pebble case, SEMFEM is added
to the considered search space.
The parameters in Table \ref{table:auto-tuner-parameter-space} are chosen as they include optimal or near-optimal preconditioner settings
for the results in Figs. \ref{fig:scaling-study-bsb}, \ref{fig:scaling-study-pb},
while still restricting the search space to something that is amenable to exhaustive search.
The authors note, however, that this parameter space may not reflect the various factors
affecting the performance of the preconditioners at especially large $P$ or on different machine architectures.

During the first timestep, our simple auto-tuner performs an exhaustive search over
the small parameter space identified in
Table \ref{table:auto-tuner-parameter-space}.
While this space is limited compared
to the true search space, the resultant preconditioners selected are effective.
In the 146 pebble case (Fig. \ref{fig:scaling-study-pb}a), Cheby-RAS(2),(7,3,1)
is identified as the preconditioner on each of the processor counts,
which was comparable in peformance to the best preconditioner.
The auto-tuner chose the optimal SEMFEM for the 1568 pebble case (Fig. \ref{fig:scaling-study-pb}b).
SEMFEM was identified at the preconditioner for the 67 pebble case on all processor counts,
which Fig. \ref{fig:scaling-study-pb}c confirms.
In the Boeing speed bump case (Fig. \ref{fig:scaling-study-bsb}), the auto-tuner
chose Cheby-RAS(2),(9,5,1) across all processor counts, which was either the optimal
or near-optimal preconditioner for the problem.

Note that the auto-tuner selects the fastest method at a {\em fixed-processor
count} (i.e., whatever the user has selected).
Consider, for example, the 67 pebble case on $P=18$ GPUs (dashed black line, Fig. \ref{fig:scaling-study-pb}c).
The goal of the auto-tuner is to select the preconditioner with lowest time-to-solution
along the user-specified iso-processor count line.
However, in the results discussed above, there was no change with respect to a change in the processor count.
We note that our principal
objective is not to squeeze out a few percent over a raft of good choices, but
rather to ensure that a case does not run with an unfortunate set of parameters
for which the performance is significantly substandard.  This primitive
auto-tuning technique proves effective at preventing the selection of highly
suboptimal preconditioners.  We have implemented this for production use and
will continue to update and refine the strategy moving forward.

\section{Conclusions and Future Work}
\label{conclusions}
In this work, we introduce Chebyshev-accelerated ASM and RAS smoothers for use
in $p$-multigrid preconditioners for the spectral element Poisson problem.
Further, we compare the performance of Schwarz-based smoothers,
Chebyshev-accelerated Jacobi smoothing, and SEMFEM-based preconditioning for a
suite of challenging test problems.  We conclude that the Chebyshev-accelerated
Schwarz smoothers with $p$-multigrid, as well as the low-order SEMFEM
preconditioner solved with AmgX, are performant on GPU machines such as OLCF's
Summit.
The authors propose a runtime auto-tuner preconditioner strategy that,
while primitive, can choose reasonable solver parameters.

The authors plan on conducting similar parallel scalability studies on other machines, such as OLCF's Spock, an
AMD MI100 machine with similar hardware and software as the upcoming Frontier system.
Additional preconditioner options, such as
varying the AMG solver
settings for the coarse grid solve as well as the solver location (CPU/GPU),
varying the Chebyshev order and number of sweeps on each level of
the multigrid heirarchy,
and varying the Chebyshev eigenvalue bounds, are avenues for future
performance optimization.
Local Fourier analysis similar to that done by Thompson {\em et al.}
in \cite{thompson2021local} is further needed to understand the
smoothing properties of the Chebyshev-accelerated Schwarz smoothers,
allowing for
robust optimization of parameters, as in the work by Brown {\em et al.} \cite{brown_tuning_2021}.

\section{Acknowledgements}
\label{acknowledgements}
This research is supported by the Exascale Computing Project (17-SC-20-SC),
a collaborative effort of two U.S. Department of Energy organizations (Office of Science and the National Nuclear
Security Administration) responsible for the planning and preparation of a capable exascale ecosystem,
including software, applications, hardware,
advanced system engineering and early testbed
platforms, in support of the nation’s exascale computing imperative.
This research also used resources of the Oak Ridge Leadership Computing Facility at Oak Ridge National
Laboratory, which is supported by the Office of
Science of the U.S. Department of Energy under
Contract DE-AC05-00OR22725.

The authors thank YuHsiang Lan, Ramesh Balakrishnan,
David Alan Reger, and Haomin Yuan
for providing visualizations and mesh files.
The authors thank the reviewers of this work for their insightful comments and suggestions.
%\vspace{1cm}
\pagebreak

%\pagebreak
\bibliographystyle{siamplain}
\bibliography{paper}

\end{document}